\documentclass[10pt]{amsart}
\usepackage{amssymb}
\usepackage{clrscode}
\usepackage{epsfig}
\usepackage{latexsym}
\usepackage{amsmath,amssymb,amsfonts,amsthm,graphics}
\usepackage{eucal}
\usepackage{eufrak}
\usepackage[all]{xypic}
\usepackage{xspace}

\theoremstyle{plain}
\newtheorem{theorem}{Theorem}

\newtheorem{lemma}{Lemma}

\theoremstyle{definition}

\theoremstyle{remark}

\numberwithin{equation}{section}

\begin{document}


\title[Asymptotic analysis of $k$-noncrossing matchings]
      {Asymptotic analysis of $k$-noncrossing matchings}
\author{Emma Y. Jin$^{}$, Christian M. Reidys$^{}$$^{\,\star}$
        and Rita R. Wang $^{}$}
\address{ Center for Combinatorics, LPMC-TJKLC \\
          Nankai University  \\
          Tianjin 300071\\
          P.R.~China\\
          Phone: *86-22-2350-6800\\
          Fax:   *86-22-2350-9272}
\email{reidys@nankai.edu.cn}
\thanks{}
\date{March, 2008}
\begin{abstract}
In this paper we study $k$-noncrossing matchings. A $k$-noncrossing
matching is a labeled graph with vertex set $\{1,\dots,2n\}$
arranged in increasing order in a horizontal line and vertex-degree
$1$. The $n$ arcs are drawn in the upper halfplane subject to the
condition that there exist no $k$ arcs that mutually intersect. We
derive: (a) for arbitrary $k$, an asymptotic approximation of the
exponential generating function of $k$-noncrossing matchings
$F_k(z)$. (b) the asymptotic formula for the number of
$k$-noncrossing matchings $f_{k}(n) \, \sim  \,c_k  \,
n^{-((k-1)^2+(k-1)/2)}\, (2(k-1))^{2n}$ for some $c_k>0$.
\end{abstract}
\maketitle {{\small
}}


\section{Statement of results and background}\label{S:intro}


Let $F_k(z)$ denote the exponential generating function of $k$-noncrossing
matchings, i.e.~
\begin{equation}
F_k(z)=\sum_{n\ge 0}f_k(n)\, \frac{z^{2n}}{(2n)!} \ .
\end{equation}
In this paper we prove the following two theorems:
\begin{theorem}\label{T:A}

Then we have for arbitrary $k\in\mathbb{N}$, $k\ge 2$, $\text{\rm
arg($z$)} \ne \pm\frac{\pi}{2}$ 
\begin{eqnarray}\label{E:das1}
F_k(z)  & = &
\left[\prod_{i=1}^{k-1}\Gamma(i+1-\frac{1}{2}) \, \prod_{r=1}^{k-2}r!\right]
\left(\frac{e^{2z}}{\pi}\right)^{k-1}\, z^{-(k-1)^2-\frac{k-1}{2}}\;
\ (1+{O}(|z|^{-1})) \ .
\end{eqnarray}
\end{theorem}
\begin{theorem}\label{T:B}
For arbitrary $k\in\mathbb{N}$, $k\ge 2$ we have
\begin{equation}\label{E:theorem}
f_{k}(n) \, \sim  \, c_k  \, n^{-((k-1)^2+(k-1)/2)}\,
(2(k-1))^{2n},\qquad \text{\it for some $c_k>0$} \ .
\end{equation}
\end{theorem}
Here we use the notation $f(z)=O(g(z))$ and $f(z)=o(g(z))$ for
$|f(z)|/|g(z)|$ being bounded and tending to zero, for $\vert z
\vert\rightarrow \infty$, respectively.

A $k$-noncrossing matching is a labeled graph over the vertices $1,\dots,2n$,
of degree exactly $1$ and drawn in increasing order in a horizontal line.
The arcs are drawn in the
upper halfplane subject to the condition that there are no $k$ arcs that
mutually intersect.
\begin{figure}[ht]
\centerline{ \epsfig{file=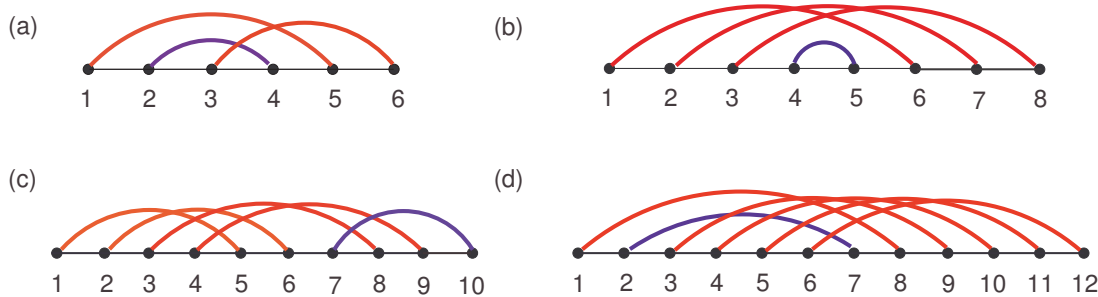,width=1.0\textwidth}\hskip15pt
} \caption{\small $k$-noncrossing matchings: $3$-, $4$-, $5$- and
$6$-noncrossing matchings respectively. One of the $k-1$ mutually
crossing arcs are drawn in red.} \label{F:diagram}
\end{figure}
Grabiner and Magyar proved an explicit determinant formula,
\cite{GM} (see also \cite{CDDSY}, eq.~$(9)$) which expresses the
exponential generating function of $f_k(n)$, for fixed $k$, as a
$(k-1) \times (k-1)$ determinant
\begin{eqnarray}
\label{E:ww1} F_k(z)=\sum_{n\ge 0} f_{k}(n)\cdot\frac{z^{2n}}{(2n)!}
& = & \det[I_{i-j}(2z)-I_{i+j}(2z)]|_{i,j=1}^{k-1} \ ,
\end{eqnarray}
where $I_m(2z)$ is the hyperbolic Bessel function:
\begin{equation}\label{E:bessel}
I_m(2z)=\sum_{j=0}^{\infty} { {z^{m+2j}} \over {j!(m+j)!} } \quad .
\end{equation}
Chen {\it et.al.} proved  in \cite{CDDSY} a beautiful correspondence
between $k$-noncrossing matchings and oscillating tableaux. The
particular RSK-insertion used in \cite{CDDSY} is based on an idea of
Stanley.
Our second result is related to a theorem of Regev \cite{Regev:81}
for the coefficient $u_k(n)$ of Gessel's generating function \cite{Gessel}
$$
U_k(x)=\det(I_{i-j}(2z))_{i,j=1}^k \ .
$$
Regev shows
\begin{equation}
u_k(n)
\sim 1!2!\dots (k-1)!\left(\frac{1}{\sqrt{2\pi}}\right)^{k-1}\,
\left(\frac{1}{2}\right)^{(k-1)^2/2} k^{k^2/2}\frac{k^{2n}}{n^{(k^2-1)/2}} \ .
\end{equation}
The proof is obtained employing the RSK-algorithm and using the
hook-length formula. One arrives, taking the limit $n\to\infty$, at
a $k$-dimensional {\it Selberg}-integral, which can be explicitely
computed. We shall use a different strategy. One key element in our
approach is the following approximation of the Bessel-function,
valid for $-\frac{\pi}{2}<
\text{\rm arg(z)} < \frac{\pi}{2}$ 
\cite{NBS:70}
\begin{equation}\label{E:bessel2}
I_m(z)=\frac{e^{z}}{\sqrt{2 \pi
z}}\left(\sum_{h=0}^{H}\frac{(-1)^{h}}{h!8^h}\prod\limits_{t=1}^h
(4m^2-(2t-1)^2)z^{-h}+\mathrm{O}(|z|^{-H-1})\right) \ .
\end{equation}
In this paper we will show that, using the approximation of
eq.~(\ref{E:bessel2}), the determinant of Bessel-functions of
eq.~(\ref{E:ww1}) can be computed asymptotically for arbitrary $k$.
The computation of the determinant via the algorithm given in
Section~\ref{S:Proof-A} is the key ingredient for all our results.


\section{Proof of Theorem~\ref{T:A}}\label{S:Proof-A}

Suppose we are given a polynomial
\begin{equation}\label{E:def-f}
g_n(x,y) = \sum_{0\leq a+b\leq n}C(a,b)\, x^{2a}y^{2b} \ ,
\end{equation}
where for $a+b=n$, $C(a,b)>0$ holds.
In the following $a$, $b$ always denote integers greater or equal to zero.
We set
\begin{equation}\label{E:triangle}
z\triangle z'= (z-z')(z+z') \ .
\end{equation}
\begin{lemma}\label{L:polynomial}
Suppose $n\ge 0$, then we have
\begin{equation}\label{E:(a)}
g_n(x,y)-g_n(x,z) =
\begin{cases}
(y\triangle z)\,
\sum\limits_{a+b\leq n-1}E(a,b,z)\, x^{2a}\,y^{2b} & n\ge 1 \\
0                                                  & n=0
\end{cases}
\end{equation}
where $E(a,b,z)=C(a,b+1)$ for $a+b=n-1$. Furthermore
\begin{equation}\label{E:(b)}
g_n(x,y)-g_n(x,y_1)-g_n(x_1,y)+g_n(x_1,y_1)=
\begin{cases}
(x\triangle x_1)(y\triangle y_1)\sum\limits_{ a+b\leq
n-2}D(a,b,x_1,y_1)x^{2a}y^{2b} & n\geq 2 \\
0 & n=0,1
\end{cases}
\end{equation}
where $D(a,b,x_1,y_1)=C(a+1,b+1)$ for $a+b=n-2$.
\end{lemma}
\begin{proof}
For $n=0$, we immediately obtain $g_0(x,y)-g_0(x,z)=0$. In case of
$n\geq 1$ we compute
\begin{align*}
g_n(x,y)-g_n(x,z) &= (y\triangle z)\sum\limits_{ a+b\leq n, \,b>0}
C(a,b)\, x^{2a}\left[\sum_{m=0}^{b-1}\, y^{2m}z^{2b-2-2m}\right]\\
&=(y\triangle z)\sum\limits_{ a+b\leq n-1}E(a,b,z)\, x^{2a}y^{2b}.
\end{align*}
In particular, for $a+b=n-1$, we observe $E(a,b,z)=C(a,b+1)$. As for
eq.~(\ref{E:(b)}) we compute in case of $n=0\ \text{or}\ 1$,
$\vartheta_n(x,x_1,y,y_1)=g_n(x,y)-g_n(x,y_1)-g_n(x_1,y)+g_n(x_1,y_1)=0$.
For $n\geq 2$ we compute
\begin{eqnarray*}
\vartheta_n(x,x_1,y,y_1)&=& (x\triangle x_1)(y\triangle
y_1)\sum\limits_{a+b\leq n, \, ab>0}
C_h(a,b)\left[\sum_{m=0}^{a-1}x^{2m}x_1^{2a-2-2m}\right]
\left[\sum_{m=0}^{b-1}y^{2m}y_1^{2b-2-2m}\right]\\
&=&(x\triangle x_1)(y\triangle y_1)\sum_{a+b\leq n-2}
D(a,b,x_1,y_1)\, x^{2a}y^{2b} \ .
\end{eqnarray*}
In particular, for $a+b=n-2$, $D(a,b,x_1,y_1)=C(a+1,b+1)$ holds.
\end{proof}
Let
\begin{eqnarray}\label{E:e}
e_{i,j}(z) & = & \sum_{h\geq 0}m_h(i,j)\frac{(-1)^h}{16^hh!}\, z^{-h}\\
\label{E:m}
m_h(i,j) & = & \prod_{t=1}^{h}(4(i-j)^2-(2t-1)^2)-
               \prod_{t=1}^{h}(4(i+j)^2-(2t-1)^2) \ .
\end{eqnarray}
We consider the algorithm {\bf A}, specified in Figure~\ref{F:algo}, which
manipulates the matrix of Laurent series
$M=(e_{i,j}(z))_{1\le i,j\le k-1}$.
\begin{figure}
\begin{codebox}
\zi $\qquad$ $\qquad$$\qquad$ $\qquad$$\qquad$ $\qquad$$\qquad$ $\qquad$
$\qquad$ $\qquad$\centerline{\bf The algorithm {\bf A}: }
\zi \kw{begin}
\zi   $M:=[e_{i,j}(z)]_{i,j=1}^{k-1};$
\zi $\quad$ $\quad$ \For $t$ from $1$ to $k-1\ \kw{do}$
\zi $\quad$$\qquad$ $\qquad$ \For $i$ from $t+1$ to $k-1$ \kw{do}
\zi $\quad$$\quad$$\qquad$ $\quad$$\qquad$ \For $j$ from $1$ to $k-1$ \kw{do}
\zi $\quad$$\quad$$\qquad$ $\quad$$\qquad$ $e_{i,j}'(z):=
\frac{-i\,
\prod_{r=1}^{t-1}(i\triangle r)}{(2t-1)!}e_{t,j}(z)+e_{i,j}(z);$
\zi $\quad$$\quad$$\qquad$ $\quad$$\qquad$ $e_{i,j}(z):=e_{i,j}'(z);$
\zi $\quad$$\quad$$\qquad$ $\quad$$\qquad$ \kw{end};
\zi $\quad$$\qquad$ $\qquad$ \kw{end};
\zi $\quad$$\qquad$ $\qquad$ \For $j$ from $t+1$ to $k-1$ \kw{do}
\zi $\quad$$\quad$$\qquad$ $\quad$$\qquad$ \For $i$ from $1$ to $k-1$ \kw{do}
\zi $\quad$$\quad$$\qquad$ $\quad$$\qquad$ $e_{i,j}''(z):=\frac{-j\,
\prod_{r=1}^{t-1}(j\triangle r)}{(2t-1)!}e_{i,t}(z)+e_{i,j}(z);$
\zi $\quad$$\quad$$\qquad$ $\quad$$\qquad$ $e_{i,j}(z):=e_{i,j}''(z);$
\zi $\quad$$\quad$$\qquad$ $\quad$$\qquad$ \kw{end};
\zi $\quad$$\qquad$ $\qquad$ \kw{end};
\zi $\quad$ $\quad$ \kw{end};
\zi  \kw{output} $M$;
\zi \kw{end};
\end{codebox}
\caption{\small }
\label{F:algo}
\end{figure}
Let $e_{i,j}^t(z)$ denote the matrix coefficient after running {\bf A} exactly
$t$ steps. We set
\begin{equation}
e_{i,j}^t(z) =  \sum_{h\geq 0} m_h^t(i,j)\frac{(-1)^h}{16^hh!}\, z^{-h}
\end{equation}
and proceed by analyzing the terms $m_h^t(i,j)$ for $1\le t<k-1$.
\begin{lemma}\label{L:induction}
For any positive integer $t$ strictly smaller than $k-1$ and we have
$m_h^t(i,j)=m_h^t(j,i)$ the following
two assertions hold.\\
{\rm (a)} for $i\leq t < j$, we have
\begin{eqnarray}
m_h^t(i,j) & = & -(2i-1)!j\,\prod_{r=1}^t(j\triangle r)\,
\sum\limits_{ a+b\leq h-(t+i)} E_h^{t}(a,b,i) \, i^{2a}j^{2b} \\
h<t+i       & \Longrightarrow & m_h^t(i,j)=0 \\
a+b=h-(t+i) & \Longrightarrow & E_h^t(a,b,i)=C_h(a+i-1,b+t)>0 \ .
\end{eqnarray}
Furthermore {\rm (a)} implies the case $j \le t<i$.

{\rm (b)} for $i,j>t$, we have
\begin{eqnarray}
m_h^t(i,j) & = & -i\,j\; \prod_{r=1}^t (i\triangle r)\,
\prod_{r=1}^t (j\triangle r)\,
\sum\limits_{ a+b\leq h-(2t+1)} D_h^{t}(a,b)\, i^{2a}j^{2b}\\
h<2t+1         & \Longrightarrow & m_h^t(i,j)=0  \\
a+b=h-(2t+1)   & \Longrightarrow & D_h^t(a,b)=C_h(a+t,b+t)>0 \ .
\end{eqnarray}
\end{lemma}
\begin{proof}
We shall prove (a) and (b) by induction on $1\le t<k-1$.
We first observe that, in view of eq.~(\ref{E:m})
\begin{equation}\label{E:qw}
\begin{split}
m_h(i,j) &=
-2\sum_{s=0}^{\lfloor\frac{h-1}{2}\rfloor}\sum_{p+q+r+2s+1=h}
\binom{h-r}{p}\binom{h-p-r}{q}(4i^2)^p(4j^2)^qa_1\cdots a_r(8ij)^{2s+1}\\
&= -ij\sum_{a+b\leq h-1}C_h(a,b)i^{2a}j^{2b} \ ,
\end{split}
\end{equation}
where $a_i\in\{-1^2,\ldots,-(2h-1)^2\}$, $i\neq j,\ a_i\neq a_j$ and
$C_h(a,b)>0$ for $a+b=h- 1$. Furthermore by definition
$m_h(i,j)=m_h(j,i)$. For $i=1,j>1$, only the $j$-loop is executed,
whence
\begin{equation}
m_h^1(i,j)= m_h(i,j)-jm_h(i,1)
\end{equation}
and for $m_h(j,i)$ only the $i$-loop contributes
$$
m_h^{1}(j,i)= m_h(j,i)-jm_h(1,i)= m_h(i,j)-jm_h(i,1)=m_h^{1}(i,j) \
.
$$
Consequently
$$
m_h^1(i,j)=-ij\left[ \sum_{a+b\leq
h-1}C_h(a,b)i^{2a}j^{2b}-\sum_{ a+b\leq h-1
}C_h(a,b)i^{2a}\right]
$$
Employing Lemma \ref{L:polynomial} we obtain
$$
m_h^1(i,j)=-i\, j\, (j\triangle 1)
\sum_{a+b\leq h-2}E_h^1(a,b,1)i^{2a}j^{2b} \ .
$$
Furthermore
\begin{eqnarray*}
a+b=h-2 & \Longrightarrow & E_h^1(a,b,1)=C_h(a,b+1)>0 \\
h= 1    & \Longrightarrow & m_h^1(i,j)=0 \ .
\end{eqnarray*}
Thus for $t=1$, the induction basis for (a) holds.
We proceed by proving that for $t=1$ (b) holds.
For $i>1,\ j>1$ both $i$- and $j$-loop are executed
\begin{eqnarray}
m_h^1(i,j)&=&m_h(i,j)-im_h(1,j)-jm_h(i,1)+ijm_h(1,1)
\end{eqnarray}
from which immediately $m_h^1(i,j)=m_h^1(j,i)$ follows. We compute
\begin{eqnarray*}
m_h^1(i,j) &=&-ij\left(\sum_{ a+b\leq h-1
}C_h(a,b)i^{2a}j^{2b}-\sum_{ a+b\leq h-1}C_h(a,b)j^{2b}\right)\\
&&-ij\left(-\sum_{ a+b\leq h-1 }C_h(a,b)i^{2a}+\sum_{
a+b\leq h-1}C_h(a,b)\right)\\
&=&-i(i\triangle 1)j(j\triangle 1)
\sum_{a+b\leq h-3 }D_h^{1}(a,b)i^{2a}j^{2b} \ .\\
\end{eqnarray*}
Lemma~\ref{L:polynomial} implies for $a+b=h-3$, $D_h^1(a,b)=C_h(a+1,b+1)$
and for $h<3$,  $m_h^1(i,j)=0$.
Accordingly we eststablished the induction basis for assertions (a),
(b) and $m_h^1(i,j)=m_h^1(j,i)$.\\
As for the induction step, we first prove {\rm (a)}. Let us suppose
{\rm (a)} holds for $t=n$. We consider the case $t=n+1$ by
distinguishing subsequent two cases: {\rm (1)} $i=n+1, j>n+1$ and
{\rm (2)} $i\le n, j=n+1$. First we observe that since $i<n+2$ the
algorithm executes no $i$-loop and by construction the only
contribution to $m_h^{n+1}(i,j)$ is made by the term
$$
-\frac{j\prod_{r=1}^{n+1-1}(j\triangle r)}{(2n+1)!} \, m_h^n(i,n+1)
$$
We accordingly derive
\begin{eqnarray*}
m_h^{n+1}(i,j) &=&
m_h^n(i,j)-\frac{j\prod_{r=1}^n(j\triangle r)}{(2n+1)!}m_h^n(i,n+1)
\end{eqnarray*}
The induction hypothesis on $t=n$ shows $m_h^n(i,j)=m_h^n(j,i)$
and $m_h^n(i,n+1)=m_h^n(n+1,i)$. Therefore we arrive at
$m_h^{n+1}(i,j)=m_h^{n+1}(j,i)$.

{\rm (1)} $i=n+1, j>n+1$, the induction hypothesis guarantees
\begin{eqnarray*}
m_h^{n}(i,j) & = & -i\,j\; \prod_{r=1}^n (i\triangle r)\,
\prod_{r=1}^{n} (j\triangle r)\,
\sum\limits_{ a+b\leq h-(2n+1)} D_h^{n}(a,b)\, i^{2a}j^{2b}\\
m_h^n(i,n+1)   & = & -i\,(n+1)\; \prod_{r=1}^n (i\triangle r)\,
\prod_{r=1}^{n} ((n+1)\triangle r)\, \sum\limits_{a+b\leq
h-(2n+1)} D_h^{n}(a,b)\, i^{2a}(n+1)^{2b}
\end{eqnarray*}
Since $(n+1)\, \prod_{r=1}^{n}((n+1)\triangle r)=(2n+1)!$
we arrive at
\begin{eqnarray*}
m_h^{n+1}(i,j)
&=&- i\,j\, \prod_{r=1}^n(i\triangle r)\, \prod_{r=1}^n(j\triangle r)
\times \\
&&
\left(\sum_{a+b\leq h-(2n+1)}D_h^n(a,b)i^{2a}j^{2b}-
\sum_{a+b\leq h-(2n+1)} D_h^n(a,b)i^{2a}(n+1)^{2b}\right)\ .\\
\end{eqnarray*}
According to Lemma \ref{L:polynomial}, for $h=2n+1$, we have
$m_h^{n+1}(n+1,j)=0$ and for $h\ge 2n+2$ we obtain
\begin{eqnarray}
\quad m_h^{n+1}(n+1,j) &=&-(2n+1)!j\prod_{r=1}^{n+1}(j\triangle r )
\sum_{a+b\leq h-(2n+2)}E_h^{n+1}(a,b,n+1)(n+1)^{2a}j^{2b}.
\end{eqnarray}
For $a+b=h-(2n+2)$, we have
\begin{equation}
E_h^{n+1}(a,b,n+1)=D_h^{n}(a,b+1)=C_h(a+n,b+n+1)>0 \ .
\end{equation}
{\rm (2)} $i\leq n$ and $j>n+1$, using the induction hypothesis, we
derive
\begin{eqnarray*}
m_h^{n+1}(i,j)&=&m_h^n(i,j)-\frac{j\prod_{r=1}^n (j\triangle
r)}{(2n+1)!}m_h^n(i,n+1)\\
 &=&-(2i-1)!\, j\, \prod_{r=1}^n(j\triangle r) \ \times \\
&& \left(\sum_{ a+b\leq h-(n+i) }E_h^{n}(a,b,i)i^{2a}j^{2b}-\sum_{
a+b\leq h-(n+i)}E_h^{n}(a,b,i)i^{2a}(n+1)^{2b}\right) \ .
\end{eqnarray*}
Lemma~\ref{L:polynomial} implies
\begin{eqnarray*}
m_h^{n+1}(i,j)&=&-(2i-1)!\, j\, \prod_{r=1}^{n+1}(j\triangle
r)\sum_{a+b\leq h-(n+1+i)}E_h^{n+1}(a,b,i)i^{2a}j^{2b}.
\end{eqnarray*}
For $h\leq n+i$, we observe $m_h^{n+1}(i,j)=0$ and for $a+b=h-(n+1+i)$,
\begin{equation}
E_h^{n+1}(a,b,i)=E_h^n(a,b+1,i)=C_h(a+i-1,b+n+1)>0 \ .
\end{equation}
Accordingly, we have proved
\begin{equation}
\forall\; i\leq n+1< j;\quad
m_h^{n+1}(i,j)=-(2i-1)!j\prod_{r=1}^{n+1}(j\triangle r)\sum_{
a+b\leq h-(n+1+i)}E_h^{n+1}(a,b,i)i^{2a}j^{2b} \ .
\end{equation}
Furthermore
\begin{eqnarray}
a+b=h-(n+1+i) & \Longrightarrow & E_h^{n+1}(a,b,i)=C_h(a+i-1,b+n+1)>0 \\
h<n+1+i & \Longrightarrow & m_h^{n+1}(i,j)=0 \ .
\end{eqnarray}
Whence assertion (a) holds by induction for any $1\le t< k-1$. We
next suppose assertion (b) is true for $t=n$ and consider the case
$t=n+1$, i.e., $i>n+1$ and $j>n+1$. First the $i$-loop is executed
and produces
$$
\tilde{m}_h^{n+1}(i,j) = m_h^n(i,j)-\frac{i\,\prod_{r=1}^n(i\triangle r)
}{(2n+1)!}m_h^n(n+1,j) \ .
$$
Secondly the $j$-loop yields
\begin{eqnarray}
m_h^{n+1}(i,j) & = & \tilde{m}_h^{n+1}(i,j)-\frac{j\prod_{r=1}^n(j\triangle r)}
{(2n+1)!}\tilde{m}_h^{n+1}(i,n+1) \ .
\end{eqnarray}
We accordingly compute
\begin{eqnarray*}
m_h^{n+1}(i,j) & = & m_h^n(i,j)-\frac{i\,\prod_{r=1}^n(i\triangle r)
}{(2n+1)!}m_h^n(n+1,j)\\
&&-\frac{j\prod_{r=1}^n(j\triangle r)}{(2n+1)!}
\left(m_h^n(i,n+1)-\frac{i\prod_{r=1}^n(i\triangle r)}
{(2n+1)!}m_h^n(n+1,n+1)\right)
\end{eqnarray*}
from which we immediately observe that $m_h^{n+1}(i,j)=m_h^{n+1}(j,i)$ holds.
Furthermore
\begin{eqnarray*}
m_h^{n+1}(i,j)&=&-i\ \left[\prod_{r=1}^n(i\triangle r)\right]\,j\,
\left[\prod_{r=1}^n(j\triangle r)\right]\ \times\\
&&\sum_{ a+b\leq h-(2n+1)}D_h^{n}(a,b)
(i^{2a}j^{2b}-(n+1)^{2a}j^{2b}-i^{2a}(n+1)^{2b}+(n+1)^{2a}(n+1)^{2b})\\
&=&-i\ \left[\prod_{r=1}^{n+1}(i\triangle r)\right]\,j\,
\left[\prod_{r=1}^{n+1}(j\triangle r)\right]
\sum_{ a+b\leq h-(2(n+1)+1)}D_h^{n+1}(a,b)i^{2a}j^{2b}.\\
\end{eqnarray*}
In particular,
\begin{eqnarray*}
a+b=h-(2(n+1)+1) &\Rightarrow &
D_h^{n+1}(a,b)=D_h^n(a+1,b+1)=C_h(a+n+1,b+n+1)>0 \\
h=2n+1 \text{ or } h=2n+2 &\Rightarrow &  m_h^{n+1}(i,j)=0 \ .
\end{eqnarray*}
Thus $m_h^{n+1}(i,j)$ satisfies (b) for any $1\le t< k-1$.
\end{proof}
We proceed by analyzing the Laurent series
\begin{equation}
a_{i,j}(z) =  \sum_{h\geq 0} m_h^{k-2}(i,j)\frac{(-1)^h}{16^hh!}\, z^{-h}\ .
\end{equation}
\begin{lemma}\label{L:asymptotics}
\begin{eqnarray}
a_{i,j}(z) & = &
(-1)^{i+j}\frac{2\Gamma(j+i-\frac{1}{2})}{\sqrt{\pi}}
z^{-(j+i-1)}(1+{O}(|z|^{-1})) \ .
\end{eqnarray}
\end{lemma}
\begin{proof} We shall prove the lemma distinguishing the cases $i<j$ and
$i=j$. The former implies by symmetry the case $i>j$. Suppose first
$i<j$. By construction of {\bf A}, we have
\begin{equation}
m_h^{k-2}(i,j)=m_h^{j-1}(i,j)
\end{equation}
since after the $(j-1)$th step, $m_h^{j-1}(i,j)$ remains unchanged. Consequently
we can write $a_{i,j}(z)$ as
\begin{equation}
a_{i,j}(z)
=\sum_{0\le h\le i+j-1}\frac{(-1)^h}{16^h\,h!}\, m_h^{j-1}(i,j)\,z^{-h}
 + \sum_{ i+j-1<h}\frac{(-1)^h}{16^h\,h!}\, m_h^{j-1}(i,j)\,z^{-h} \ .
\end{equation}
We consider the terms $m_h^{j-1}(i,j)$ for
$0\le h\le j+i-1$. According to Lemma~\ref{L:induction}
\begin{equation*}
m_h^{j-1}(i,j)=-(2i-1)!\, j\, \prod_{r=1}^{j-1}(j\triangle r)\, \sum_{0\leq
a+b\leq h-(j-1+i)}E_h^{j-1}(a,b,i)\, i^{2a}j^{2b}
\end{equation*}
holds. In particular,
\begin{equation}
h<j-1+i \quad \Longrightarrow \quad m_h^{j-1}(i,j)=0 \ .
\end{equation}
Accordingly, the only nonzero coefficient of
$\sum_{0\le h\le i+j-1}\frac{(-1)^h}{16^h\,h!}\, m_h^{j-1}(i,j)\,z^{-h}$
has index $h=j-1+i$ in which case
$$
a+b=0\quad \text{\rm and}\quad  E_{j-1+i}^{j-1}(0,0,i)=C_{j-1+i}(i-1,j-1)
$$
holds, i.e.
\begin{equation}
a_{i,j}(z)= \frac{(-1)^{j+i}C_{j-1+i}(i-1,j-1)(2j-1)!(2i-1)!}
{16^{j-1+i}(j-1+i)!}z^{-(j-1+i)}\, (1+{O}(|z|^{-1})) \ .
\end{equation}
Secondly suppose $i=j$. Then, by definition of ${\bf A}$, the Laurent
series $a_{i,i}(z)$ is obtained for $t=i-1$, i.e.~we have
\begin{equation}
a_{i,i}(z)
=\sum_{0\le h\le 2i-1}\frac{(-1)^h}{16^h\,h!}\, m_h^{i-1}(i,i)\,z^{-h}
 + \sum_{ 2i-1<h}\frac{(-1)^h}{16^h\,h!}\, m_h^{i-1}(i,i)\,z^{-h} \ .
\end{equation}
Lemma~\ref{L:induction} (b) implies
\begin{equation*}
m_h^{i-1}(i,i)=-((2i-1)!)^2\sum_{0\leq a+b \leq
h-(2i-1)}D_h^{i-1}(a,b)\, i^{2a}i^{2b}.
\end{equation*}
In particular for $h<2i-1$ we have $m_h^{i-1}(i,i)=0$, thus for
$\sum_{0\le h\le 2i-1}\frac{(-1)^h}{16^h\,h!}\, m_h^{i-1}(i,i)\,z^{-h}$ only
the index $h=2i-1$ has a nonzero coefficient in which case
$$
a=b=0 \quad \text{\rm and} \quad D_{2i-1}^{i-1}(0,0)=C_{2i-1}(i-1,i-1)
$$
holds. We therefore derive
\begin{equation}
a_{i,i}(z)= \frac{(-1)^{2i}((2i-1)!)^2C_{2i-1}(i-1,i-1)}
{16^{2i-1}(2i-1)!}z^{-(2i-1)}(1+{O}(|z|^{-1})).
\end{equation}
Thus we have proved that we have for $i\leq j$
\begin{equation}
a_{i,j}(z)= \frac{(-1)^{j+i}C_{j-1+i}(i-1,j-1)(2j-1)!(2i-1)!}
{16^{j-1+i}(j-1+i)!}z^{-(j-1+i)}(1+{O}(|z|^{-1})) \ .
\end{equation}
{\it Claim $1$.}
\begin{equation}
C_{j-1+i}(i-1,j-1)=
\frac{j\Gamma(2i+2j-1)}{\Gamma(2j+1)\Gamma(2i)} 4^{j+i}.
\end{equation}
According to eq.~(\ref{E:qw})
\begin{eqnarray*}
m_h(i,j)&=& -ij\sum_{0\leq a+b\leq h-1}C_h(a,b)i^{2a}j^{2b}
\end{eqnarray*}
from which we can conclude for $a+b=h-1$ and $l_m=\min\{a,b\}$
\begin{align*}
C_h(a,b)&=
\sum_{s=0}^{l_m}\binom{h}{a-s}\binom{h-a+s}{b-s}4^{a-s}4^{b-s}
8^{2s+1}2\nonumber\\
&=4^{h+1}\sum_{s=0}^{l_m}\binom{h}{a-s}\binom{h-a+s}{2s+1}4^s \ .
\end{align*}
Therefore
\begin{equation}
C_{j-1+i}(i-1,j-1)=
\frac{j\Gamma(2i+2j-1)}{\Gamma(2j+1)\Gamma(2i)} 4^{j+i} \
\end{equation}
and Claim $1$ follows. Since $\det[a_{i,j}(z)]_{i,j=1}^{k-1}$ is
symmetric, we arrive at
\begin{equation*}
a_{i,j}(z)=(-1)^{i+j}\frac{2\Gamma(j+i-\frac{1}{2})}
{\sqrt{\pi}}z^{-(j+i-1)}(1+{O}(|z|^{-1}))
\end{equation*}
for any $1\le i,j\le k-1$ and the lemma follows.
\end{proof}


{\bf Proof of Theorem~\ref{T:A}.}\\
Let
\begin{equation}
b_{i,j}(z)=(-1)^{i+j}\frac{2\Gamma(j+i-\frac{1}{2})}{\sqrt{\pi}}z^{-(j+i-1)}
\ .
\end{equation}
According to Lemma~\ref{L:asymptotics} we have $a_{i,j}(z)=
b_{i,j}(z)\, [1+{O}(|z|^{-1})]$ and we immediately obtain
\begin{eqnarray*}
\det[a_{i,j}(z)]_{i,j=1}^{k-1} &=&\sum_{\sigma\in S_{k-1}}{\rm
sign}(\sigma)\, \prod_{j=1}^{k-1}\left[
b_{j,\sigma(j)}(z)\,\left[1+{O}(|z|^{-1})\right]\right]
\end{eqnarray*}
where $S_{k-1}$ denotes the symmetric group over $k-1$ letters. Furthermore
we observe
\begin{eqnarray*}
\det[b_{i,j}(z)]_{i,j=1}^{k-1}&=&\sum_{\sigma\in S_{k-1}}
\text{\rm sign}(\sigma)\,\prod_{j=1}^{k-1}b_{j,\sigma(j)}(z)\\
&=&\sum_{\sigma\in S_{k-1}}
\text{\rm sign}(\sigma)\,(-1)^{\sum_{j=1}^{k-1}(j+\sigma(j))}\,
\left(\frac{2}{\sqrt{\pi}}\right)^{k-1} \times \\
& & \qquad \qquad z^{-\sum_{j=1}^{k-1}(j+\sigma(j)-1)}
                            \prod_{j=1}^{k-1}\Gamma(j+\sigma(j)-\frac{1}{2})\\
\end{eqnarray*}
Since $\sum_{j=1}^{k-1}(j+\sigma(j))=k(k-1)$ we arrive at
\begin{eqnarray*}
\det[b_{i,j}(z)]_{i,j=1}^{k-1}
& = & \left(\frac{2}{\sqrt{\pi}}\right)^{k-1} z^{-(k-1)^2}\,
\det\left[\Gamma(j+i-\frac{1}{2})\right]_{i,j=1}^{k-1} \\
\end{eqnarray*}
and consequently
\begin{eqnarray*}
\det[a_{i,j}(z)]_{i,j=1}^{k-1} & = & \sum_{\sigma\in S_{k-1}}{\rm
sign}(\sigma)\, \prod_{j=1}^{k-1}\left[
b_{j,\sigma(j)}(z)\,\left[1+{O}(|z|^{-1})\right]\right]\\
 & =
& \det[\Gamma(j+i-\frac{1}{2})]_{i,j=1}^{k-1}
\left(\frac{2}{\sqrt{\pi}}\right)^{k-1}
z^{-(k-1)^2}(1+{O}(|z|^{-1})) \ .
\end{eqnarray*}
We proceed by computing the determinant
\begin{equation}
\det[\Gamma(j+i-\frac{1}{2})]_{i,j=1}^{k-1} =
\prod_{i=1}^{k-1}\Gamma(i+1-\frac{1}{2}) \, \prod_{r=1}^{k-2}r! \ .
\end{equation}
Since $\Gamma(i+j+1-1/2)=(i+j-1/2)\,\Gamma(i+j-1/2)$, we have for $j>1$
\begin{equation*}
\Gamma(i+j-\frac{1}{2})=\prod_{r=1}^{j-1}(i+r-\frac{1}{2}) \,
\Gamma(i+1-\frac{1}{2}) \ .
\end{equation*}
We set
\begin{equation*}
 u_{i,j}=\begin{cases} \prod\limits_{r=1}^{j-1}(i+r-\frac{1}{2})& j>1\\
 1& j=1\\
 \end{cases}
\end{equation*}
and compute
\begin{align*}
\det[\Gamma(j+i-\frac{1}{2})]_{i,j=1}^{k-1}&=
\prod_{i=1}^{k-1}\Gamma(i+1-\frac{1}{2})
\det[u_{i,j}]_{i,j=1}^{k-1}
=\prod_{i=1}^{k-1}\Gamma(i+1-\frac{1}{2})\det[i^{j-1}]_{i,j=1}^{k-1}.
\end{align*}
The determinant $\det[i^{j-1}]_{i,j=1}^{k-1}$ is a Vandermonde
determinant, whence
\begin{equation*}
\det[i^{j-1}]_{i,j=1}^{k-1}=\sum_{1\leq i_1<i_2\leq
k-1}(i_2-i_1)=\prod_{r=1}^{k-2}r! \ .
\end{equation*}
Therefore we have shown
\begin{eqnarray}\label{E:www}
\det[a_{i,j}(z)]_{i,j=1}^{k-1} & = &
\left[\prod_{i=1}^{k-1}\Gamma(i+1-\frac{1}{2}) \, \prod_{r=1}^{k-2}r!\right]\,
\left(\frac{2}{\sqrt{\pi}}\right)^{k-1}
z^{-(k-1)^2}(1+{O}(|z|^{-1})) \ .
\end{eqnarray}
It remains to combine our findings:
the approximation of the Bessel
function eq.~(\ref{E:bessel2}) and eq.~(\ref{E:m}) imply for
$-\frac{\pi}{2}< \text{\rm arg(z)} < \frac{\pi}{2}$
\begin{equation*}
I_{i-j}(2z)-I_{i+j}(2z)=\frac{e^{2z}}{2\sqrt{\pi
z}}\left(\sum_{h=1}^H m_h(i,j)\frac{(-1)^h}{16^hh!}\,
z^{-h}+O(|z|^{-H-1})\right) \ .
\end{equation*}
Let
\begin{equation}
e_{i,j}^H(z)=\sum_{h=1}^H m_h(i,j)\frac{(-1)^h}{16^hh!}\,z^{-h}
\end{equation}
then we have
\begin{equation}
F_k(z)=\det[I_{i-j}(2z)-I_{i+j}(2z)]_{i,j=1}^{k-1}
=\left(\frac{e^{2z}}{2\sqrt{\pi z}}\right)^{k-1}\,
\left[\det\left[e_{i,j}^H\right]_{i,j=1}^{k-1}+O(|z|^{-H-1})\right].
\end{equation}
Lemma~\ref{L:induction} and Lemma~\ref{L:asymptotics} provide an
interpretation of $\det[e_{i,j}^H(z)]_{i,j=1}^{k-1}$: for
\begin{equation}
H> (k-1)^2
\end{equation}
we can conclude
\begin{equation*}
\det[e_{i,j}^H(z)]_{i,j=1}^{k-1}=\det[b_{i,j}(z)]_{i,j=1}^{k-1} \,
\left[1+O(\vert z\vert^{-1})\right] \ .
\end{equation*}
Accordingly we derive
\begin{eqnarray*}
F_k(z)  & = & \left(\frac{e^{2z}}{2\sqrt{\pi z}}\right)^{k-1}\,
\det[b_{i,j}(z)]_{i,j=1}^{k-1} \, \left[1+O(\vert
z\vert^{-1})\right] \ .
\end{eqnarray*}
Since
$$
\det[b_{i,j}(z)]_{i,j=1}^{k-1}  =
\left(\frac{2}{\sqrt{\pi}}\right)^{k-1}
z^{-(k-1)^2}\,\det[\Gamma(j+i-\frac{1}{2})]_{i,j=1}^{k-1}
$$
and $F_k(z)$ is an even function, we obtain for $\text{\rm arg}(z)\ne
\pm \frac{\pi}{2}$
\begin{eqnarray}\label{E:das}
F_k(z)  & = &
\left[\prod_{i=1}^{k-1}\Gamma(i+1-\frac{1}{2}) \, \prod_{r=1}^{k-2}r!\right]
\left(\frac{e^{2z}}{\pi}\right)^{k-1}\, z^{-(k-1)^2-\frac{k-1}{2}}\;
\ (1+{O}(|z|^{-1}))
\end{eqnarray}
and the proof of the theorem is complete. \hfill $\square$

\section{Proof of Theorem~\ref{T:B}}


Suppose $k=4m$, $m\in\mathbb{N}$,
$p=(k-1)^2+\frac{k-2}{2}=(4m-1)^2+2m-1$ and
\begin{eqnarray} \label{E:case1}
\quad g_k(z) & = & \tilde{c}_{k}\,
\left[I_0((2k-2)z)\,z^{-p}-\sum_{j=1}^{p}a_{k,j}\,z^{-j}\right], \
\text{\rm where} \ a_{k,j}=[z^{p-j}]\,I_0((2k-2)z).
\end{eqnarray}
For $k=4m+2$, let $p=(k-1)^2+\frac{k-2}{2}=(4m+1)^2+2m$ and
\begin{eqnarray} \label{E:case3}
\quad g_k(z) &= &
\tilde{c}_{k}\left[I_1((2k-2)z)\,z^{-p}-\sum_{j=1}^{p}a_{k,j}\,z^{-j}\right],
\ \text{\rm where} \ a_{k,j}=[z^{p-j}]\,I_1((2k-2)z).
\end{eqnarray}
For $k=4m+1$, let $p=(k-1)^2+\frac{k-1}{2}=(4m)^2+2m$ we set
\begin{eqnarray} \label{E:case2}
\quad g_k(z) & = & \tilde{c}_{k}\,
\left[\cosh((2k-2)z)\,z^{-p}-\sum_{j=1}^{p}a_{k,j}\,z^{-j}\right],\
\text{\rm where} \ a_{k,j}=[z^{p-j}]\,\cosh((2k-2)z).
\end{eqnarray}
Finally, for $k=4m+3$, let $p=(k-1)^2+\frac{k-1}{2}=(4m+2)^2+2m+1$
and
\begin{eqnarray} \label{E:case4}
\quad g_k(z) & = &
\tilde{c}_{k}\left[\sinh((2k-2)z)\,z^{-p}-\sum_{j=1}^{p}a_{k,j}\,z^{-j}\right],
\ \text{\rm where} \ \ a_{k,j}=[z^{p-j}]\,\sinh((2k-2)z).
\end{eqnarray}
The functions given in eq.~(\ref{E:case1})-(\ref{E:case4}) are
entire, even and the constants $\tilde{c}_{k}$ satisfy
\begin{equation*} g_k(|z|)\sim
c_k'\,e^{(2k-2)|z|}\,|z|^{-(k-1)^2-\frac{k-1}{2}},\quad \ \text{\rm
as} \ \vert z\vert\to \infty \
\end{equation*}
where $c_k'=\pi^{-(k-1)}\prod_{i=1}^{k-1}\Gamma(i+1-\frac{1}{2}) \,
\prod_{r=1}^{k-2}r!$.
\begin{proof}
{\it Claim $1$.} Suppose $z\in\mathbb{C}\setminus\mathbb{R}$, then
we have
\begin{equation}\label{E:F_kpro}
|F_k(z)|={o}(|z|^{-1}F_k(|z|)).
\end{equation}
To prove Claim $1$, we conclude from Theorem~\ref{T:A} that
\begin{equation}\label{E:F_kasy}
F_k(z)=
c_k'\,e^{(2k-2)z}\,z^{-(k-1)^2-\frac{k-1}{2}}(1+{O}(|z|^{-1}))\quad\text{\rm
for}\ \arg(z)\neq \pm \pi/2\ ,
\end{equation}
where $c_k'=\pi^{-(k-1)}\prod_{i=1}^{k-1}\Gamma(i+1-\frac{1}{2}) \,
\prod_{r=1}^{k-2}r!$ holds. We write $z=re^{i\theta}$ and obtain for
$\theta\neq 0,\ \pi,\ \pm \pi/2$
\begin{eqnarray}\label{E:fast}
\frac{|F_k(z)|}{|z|^{-1}F_k(|z|)} & = & e^{-2(k-1)(1-\cos\theta)\,
r}\ r \ \left({O}(1)+{O}(|z|^{-1})\right).
\end{eqnarray}
Therefore we have $|F_k(z)|={o}( |z|^{-1}F_k(|z|))$ for
$\arg(z)\neq0,\pi, \pm \pi/2$. Since $|F_k(z)|$ and $|z|^{-1}F_k(|z|)$
are continuous, eq.~(\ref{E:fast}) implies
\begin{eqnarray}
|F_k(z)| & = & {o}(|z|^{-1}F_k(|z|)), \quad \text{\rm for
$z\in\mathbb{C} \setminus\mathbb{R}$}.
\end{eqnarray}
whence Claim $1$.\\
{\it Claim 2.} For any $k\ge 2$, the functions given in
eq.~(\ref{E:case1})-(\ref{E:case4}) satisfy
\begin{equation}\label{E:g_kpro}
|g_k(z)|={o}(|z|^{-1}g_k(|z|)) \quad \text{\rm for $z\in\mathbb{C}
\setminus\mathbb{R}$}
\end{equation}
and
\begin{equation*} g_k(|z|)=
c_k'\,e^{(2k-2)|z|}|z|^{-(k-1)^2-\frac{k-1}{2}}(1+{O}(|z|^{-1})).
\end{equation*}
Suppose first $k=4m$ or $4m+2$. Then we have
\begin{equation}\label{E:g_keven}
g_k(z)=\tilde{c}_{k}
\left(I_s((2k-2)z)z^{-p}-\sum_{j=1}^{p}a_{k,j}z^{-j}\right),\quad
s=0\ \text{or}\ 1,
\end{equation}
where $p=(k-1)^2+\frac{k-2}{2}$. For
$-\frac{\pi}{2}<\text{arg}(z)<\frac{\pi}{2}$, we have
\begin{equation}\label{E:bessel3}
I_s(z)=\frac{e^{z}}{\sqrt{2 \pi
z}}\left(\sum_{h=0}^{H}\frac{(-1)^{h}}{h!8^h}\prod\limits_{t=1}^h
(4s^2-(2t-1)^2)z^{-h}+{O}(|z|^{-H-1})\right).
\end{equation}
Using eq.(\ref{E:bessel3}) we derive for sufficiently large $\vert
z\vert$
\begin{align*}
\frac{|g_k(z)|}{|z|^{-1}g_k(|z|)}&\leq
\frac{|I_s((2k-2)z)||z|^{-p}+\sum_{j=1}^{p}a_{k,j}|z|^{-j}}
{I_s((2k-2)|z|)|z|^{-p-1}-\sum_{j=1}^{p}a_{k,j}|z|^{-j-1}}  \\
 & \le  C_0\; e^{-2(k-1)(1-\cos\theta)r}\; r\; \ ,
\end{align*}
where $C_0>0$ is some constant. Since $g_k(z)$ is even we have shown

\begin{equation}\label{E:fast2}
|g_k(z)| ={o}(|z|^{-1}g_k(|z|)) \quad  \; \text{\rm where
}\arg(z)\not\in\{0,\pi,\frac{\pi}{2},-\frac{\pi}{2} \}.
\end{equation}
Since $g_k(z)$ is continuous eq.~(\ref{E:fast2}) implies
$|g_k(z)|={o}(|z|^{-1}g_k(|z|)) \quad \text{\rm for $z\in\mathbb{C}
\setminus\mathbb{R}$}$ .\\
By eq.~(\ref{E:bessel3}) and the definition of $g_k(z)$,  we can
obtain that
\begin{align*}
g_k(|z|)&=\tilde{c}_{k}
\left(I_s((2k-2)|z|)|z|^{-p}-\sum_{j=1}^{p}a_{k,j}|z|^{-j}\right)\\
&=\tilde{c}_{k}\frac{e^{(2k-2)|z|}}{2\sqrt{(k-1)\pi}\,|z|^{p+\frac{1}{2}}}
(1+O(|z|^{-1}))-\tilde{c}_{k}\sum_{j=1}^{p}a_{k,j}|z|^{-j}\\
&= c_k'\,e^{(2k-2)|z|}|z|^{-(k-1)^2-\frac{k-1}{2}}(1+{O}(|z|^{-1})).
\end{align*}

For $k=4m+1$ or $4m+3$, $g_k(z)$ satisfies
\begin{align*}
|g_k(z)|&\leq \tilde{c}_k\left(\frac{|e^{(2k-2)z}|+|e^{-(2k-2)z}|}
{2}|z|^{-p}+\sum_{j=1}^{p}a_{k,j}|z|^{-j}\right)\\
&=\tilde{c}_k\left(\frac{e^{(2k-2)r\cos\theta}+
e^{-(2k-2)r\cos\theta}}{2}r^{-p}+\sum_{j=1}^{p}a_{k,j}r^{-j}\right)
\end{align*}
where $p=(k-1)^2+\frac{k-1}{2}$ and consequently for sufficiently
large $\vert z\vert$
\begin{equation}\label{E:oha}
\frac{|g_k(z)|}{|z|^{-1}g_k(|z|)} \le C_1\; r\;
e^{-(2k-2)r(1-|\cos\theta|)}
\end{equation}
for some $C_1>0$. eq.~(\ref{E:oha}) shows
\begin{equation}
\forall z\in\mathbb{C}\setminus \mathbb{R} \qquad
|g_k(z)|={o}(|z|^{-1}g_k(|z|)) \ .
\end{equation}
For $k=4m+1$ we derive
\begin{align*}
g_k(|z|)&=\tilde{c}_k\left(\cosh((2k-2)|z|)|z|^{-p}
-\sum_{j=1}^{p}a_{k,j}|z|^{-j}\right)\\
&=\tilde{c}_k\left(\frac{e^{(2k-2)|z|}+e^{-(2k-2)|z|}}{2}|z|^{-p}
-\sum_{j=1}^{p}a_{k,j}|z|^{-j}\right)\\
&=c_k'\,e^{(2k-2)|z|}|z|^{-(k-1)^2-\frac{k-1}{2}}(1+{O}(|z|^{-1})).
\end{align*}
The case $k=4m+3$ follows analogously.
We can conclude from
\begin{equation*}
F_k(z) =
c_k'\,e^{(2k-2)z}z^{-(k-1)^2-\frac{k-1}{2}}(1+{O}(|z|^{-1}))\quad\text{\rm
for}\ \arg(z)\neq \pm \pi/2, \\
\end{equation*}
and
\begin{equation*} g_k(|z|)=
c_k'\,e^{(2k-2)|z|}|z|^{-(k-1)^2-\frac{k-1}{2}}(1+{O}(|z|^{-1})).
\end{equation*}
that $F_k(\vert z\vert) =  g_k(\vert z\vert)(1+{O}(|z|^{-1}))$
holds.
To summarize we have shown
\begin{eqnarray*}
|F_k(z)| & = & {o}(|z|^{-1}F_k(|z|)) \quad\text{\rm for
$z\in\mathbb{C}\setminus \mathbb{R}$}\\
|g_k(z)| & = & {o}(|z|^{-1}g_k(|z|))  \quad\text{\rm for
$z\in\mathbb{C}\setminus \mathbb{R}$}\\
F_k(\vert z\vert)& = & g_k(\vert z\vert)(1+{O}(|z|^{-1})).
\end{eqnarray*}
We can accordingly conclude that
\begin{equation}\label{E:wichtig}
|F_k(z)-g_k(z)|={O}(|z|^{-1}g_k(|z|)),
\end{equation}
uniformly for all $z$ with $|z|\geq 1$.\\
{\it Claim 3.} For arbitrary $k\ge 2$ we have
\begin{equation}\label{E:coeff33}
f_k(n)\sim c_k\, n^{-(k-1)^2-\frac{k-1}{2}}\;(2k-2)^{2n}
\quad\text{\rm where $c_k>0$} \ .
\end{equation}
To prove Claim $3$ we compute, using eq.~(\ref{E:wichtig})
\begin{eqnarray*}
\left| [z^{2n}]\,(F_k(z)-g_k(z))\right| &\leq&\int_{\vert
z\vert=\frac{n}{k-1}}\frac{\left|F_k(z)-g_k(z)\right|}{|z|^{2n+1}}|\mathrm{d}z|\\
&\leq& c\int_{\vert
z\vert=\frac{n}{k-1}}\frac{|z|^{-1}g_k(|z|)}{|z|^{2n+1}}|\mathrm{d}z|,\\
\end{eqnarray*}
where $c$ is a positive constant. For $k=4m\ \text{or}\ 4m+2$ we
have $p=(k-1)^2+\frac{k-2}{2}$ and substituting for $g_k(\vert
z\vert)$
\begin{align*}
\left| [z^{2n}]\,(F_k(z)-g_k(z))\right| &\leq c'\int_{\vert
z\vert=\frac{n}{k-1}}\frac{|z|^{-1}|z|^{-p-\frac{1}{2}}
e^{(2k-2)|z|}}{|z|^{2n+1}}|\mathrm{d}z|\\
&=c'\, e^{(2k-2)\cdot\frac{n}{k-1}}\, \left(\frac{n}{k-1}\right)^{
-2n-2-p-\frac{1}{2}}\, 2\pi\frac{n}{k-1}\\
&=c''\, e^{2n}(k-1)^{2n}n^{-(2n+p+\frac{3}{2})}
\end{align*}
where $c',c''$ are positive constants. By definition of the Bessel
function, see eq.~(\ref{E:bessel}),~(\ref{E:g_keven}) and using
Stirling's formula
\begin{equation}\label{stirling1}
\begin{split}
[z^{2n}]\,g_k(z)&=\tilde{c}_{k}\,[z^{2n+p}]\,I_s((2k-2)z)
=\tilde{c}_{k}\,\frac{(k-1)^{2n+p}}{(n+\frac{p-s}{2})!(n+\frac{p+s}{2})!}\\
&\sim \tilde{c}_{k}'\; e^{2n}(k-1)^{2n}n^{-(2n+p+1)}.
\end{split}
\end{equation}
Here $s$ only depends on $k$ and $\tilde{c}_{k}'$ is a positive
constant. Therefore we can conclude
\begin{equation}
[z^{2n}]\,F_k(z) \sim [z^{2n}]\, g_k(z) \ ,
\end{equation}
whence
\begin{eqnarray*}
f_{k}(n)  &=& (2n)!\,[z^{2n}]\, F_k(z) \, \sim\,
(2n)!\tilde{c}_{k}\,
                      \frac{(k-1)^{2n+p}}{(n+\frac{p-s}{2})!(n+\frac{p+s}{2})!}\\
 &\sim & c_k\,(2k-2)^{2n}n^{-(k-1)^2-\frac{k-1}{2}} \ .
\end{eqnarray*}
In case of $k=4m+1$ or $4m+3$ we have $p=(k-1)^2+\frac{k-1}{2}$ and
compute
\begin{align*}
\left| [z^{2n}](F_k(z)-g_k(z))\right| & \leq c'\int_{\vert
z\vert=\frac{n}{k-1}}
\frac{|z|^{-1}|z|^{-p}e^{(2k-2)|z|}}{|z|^{2n+1}}|\mathrm{d}z|\\
&=c'e^{(2k-2)\, \frac{n}{k-1}}\,
\left(\frac{n}{k-1}\right)^{-2n-2-p}
\, 2\pi\frac{n}{k-1}\\
&=c''e^{2n}(k-1)^{2n}n^{-(2n+p+1)}
\end{align*}
where $c',c''$ are positive constants. For $k=4m+1$ we obtain
\begin{equation}\label{stirling2}
[z^{2n}]\,g_k(z) =  \tilde{c}_{k} \,
[z^{2n+p}]\,\cosh((2k-2)z)= \tilde{c}_{k}\,
\frac{(2k-2)^{2n+p}}{(2n+p)!} \sim \tilde{c}_{k}'\;
e^{2n}(k-1)^{2n}n^{-(2n+p+\frac{1}{2})}
\end{equation}
and for $k=4m+3$
\begin{equation}\label{stirling3}
[z^{2n}]\,g_k(z)=\tilde{c}_{k}\,[z^{2n+p}]\,\sinh((2k-2)z)
=\tilde{c}_{k}\, \frac{(2k-2)^{2n+p}}{(2n+p)!} \sim \tilde{c}_{k}'\;
e^{2n}(k-1)^{2n}n^{-(2n+p+\frac{1}{2})}.
\end{equation}
Since $\left| [z^{2n}]\,(F_k(z)-g_k(z))\right|\le
c''\,e^{2n}(k-1)^{2n}n^{-(2n+p+1)}$ eq.~(\ref{stirling2}) and
(\ref{stirling3}) guarantee
\begin{equation}
[z^{2n}]\,F_k(z) \sim [z^{2n}]\,g_k(z) \ .
\end{equation}
Accordingly we obtain
\begin{equation}
f_{k}(n)=(2n)!\,[z^{2n}]\,F_k(z)\sim (2n)!\,\tilde{c}_{k}\,
\frac{(2k-2)^{2n+p}}{(2n+p)!}\sim c_k\, n^{-(k-1)^2-\frac{k-1}{2}}
\, (2k-2)^{2n}
\end{equation}
and Theorem~\ref{T:B} follows.
\end{proof}

{\bf Acknowledgments.}
This work was supported by the 973 Project, the PCSIRT Project of the
Ministry of Education, the Ministry of Science and Technology, and
the National Science Foundation of China.

\bibliographystyle{amsplain}


\end{document}